%
%
%
%
%
%
%
%

\input math.macros
\input Ref.macros

\newif\ifsubmit

\submitfalse
\ifsubmit
\hsize=138truemm\vsize=215truemm
\parskip=0pt\parindent=16truept
\fi

\checkdefinedreferencetrue
\continuousfigurenumberingtrue
\theoremcountingtrue
\sectionnumberstrue
\forwardreferencetrue
\initialeqmacro

\def\ssection#1{\medbreak\noindent{\bf #1}.}

\def\vertex{V}
\def\edge{E}
\def\Aut{{\rm Aut}}
\def\MTP/{Mass-Transport Principle}
\def\\{\backslash}
\def\bd{\partial}
\def\iso{\iota}
\def\bde{\bd_\edge}

\def\isoe{\iso_\edge}

\def\st{:\,}
\def\path{{\cal P}}
\def\H{{\Bbb H}}

\def\OHD{{\cal O}_{\bf HD}}

\font\frak=eufm10   %
\def\fo{{\hbox{\frak F}}}
\def\Eow{\E^\omega_\wsf}
\def\Eof{\E^\omega_\fsf}
\def\EGw{\E^G_\wsf}
\def\EGf{\E^G_\fsf}

\def\Pof{\P^\omega_\fsf}

\def\PGf{\P^G_\fsf}
\def\fsf{{\ss FSF}}
\def\wsf{{\ss WSF}}
\def\owsf{{\ss OWSF}}         %
\def\fusf{{\ss FSF}}          %
\def\wusf{{\ss WSF}}          %
\def\verts{\vertex}
\def\edges{\edge}
\def\dist{{\rm dist}}
\def\distsub{{\rm dist}}
\def\ev#1{{\cal #1}}
\def\SS{{\cal S}}   %
\def\setminu{-}
\def\le{\leqslant}
\def\leq{\leqslant}
\def\ge{\geqslant}
\def\geq{\geqslant}

\def\assumptions{Let $G$ be a graph with a transitive unimodular
                 closed automorphism group $\Gamma\subset\Aut(G)$}
\def\gr{{\rm gr}}
\def\ins{\Pi}   %
\def\anc{\isoe^*}   %
\def\itemrm#1{\item{{\rm #1}}}
\def\beginitems{\begingroup\parindent=25pt}
\def\enditems{\vskip0pt\endgroup\noindent}

\def\BLPSgip{Benjamini, Lyons, Peres, and Schramm (1997)%
\def\BLPSgip{BLPS (1997)}}
\def\BLPSusf{Benjamini, Lyons, Peres, and Schramm (1998)%
\def\BLPSusf{BLPS (1998)}}

\ifproofmode \relax \else\head{} {Version of 17 Mar.\ 1998}\fi
\vglue20pt

\title{Percolation Perturbations in}
\title{Potential Theory and Random Walks}

\author{Itai Benjamini, Russell Lyons, and Oded Schramm}

\abstract{%
We show that on a Cayley graph of a nonamenable group, a.s.\ the
infinite clusters of Bernoulli percolation are transient for
simple random walk, that simple
random walk on these clusters has positive speed, and that these
clusters admit bounded harmonic functions.
A principal new finding on which these results are
based is that such clusters admit invariant
random subgraphs with positive isoperimetric constant.
\par
We also show that percolation clusters 
in any amenable Cayley graph a.s.\ admit no nonconstant harmonic 
Dirichlet functions.
Conversely, on a Cayley graph admitting nonconstant harmonic 
Dirichlet functions, 
a.s.\ the infinite clusters of $p$-Bernoulli percolation 
also have nonconstant harmonic Dirichlet functions
when $p$ is sufficiently close to $1$.
Many conjectures and questions are presented.}

\bottom{Primary 60B99. %
 Secondary
60D05, %
31C20, %
60J15, %
31B05, %
20F32%
.} 
{Spanning trees, Cayley graphs, harmonic Dirichlet
functions, bounded harmonic functions,
amenability, percolation, random walks, speed, entropy,
RWRE (random walk in a random environment), isoperimetric constant.}
{Research supported by a Varon Visiting Professorship at
the Weizmann Institute of Science (Lyons) and
the Sam and Ayala Zacks Professorial Chair (Schramm).}

\bsection {Introduction}{s.intro}

The question of whether various potential-theoretic properties of graphs
and manifolds are preserved under perturbations or approximations has been
studied for more than a decade.
For example, invariance under quasi-isometries of transience
(i.e., the existence of nonconstant positive superharmonic functions) 
or of existence of harmonic functions in certain classes has been studied by
Kanai (1986), T. Lyons (1987), Saloff-Coste (1992), Soardi (1993),
Benjamini and Schramm (1996a), Thm.~3.5, and Holopainen and Soardi (1997).

In this paper, 
we study perturbations of graphs that are more radical than
quasi-isometries and that are random.
Namely, edges are deleted at random to form a
percolation subgraph $\omega$ and the behavior of {\bf simple
random walk} $\Seq{X(t)}$ on $\omega$ is examined
(where each neighbor of $X(t)$ in $\omega$ is equally likely to be $X(t+1)$).

We recall some definitions.
Given a graph
$G = \bigl(\verts(G), \edges(G)\bigr)$
and $p \in [0, 1]$, the random subgraph $\omega_p$
formed by deleting each edge independently with probability $1-p$ is called
{\bf $p$-Bernoulli bond percolation}.
The {\bf critical probability} $p_c(G)$ is the infimum over all $p\in[0,1]$ such
that there is positive probability
for the existence of an infinite connected component in $\omega_p$.
The connected components of $\omega_p$ are also called {\bf clusters}.
In the case that $G$ is an amenable Cayley graph, Burton and Keane (1989)
show that $p$-Bernoulli percolation has a.s.\ at most
one infinite cluster. For background on percolation, especially in
$\Z^d$, see Grimmett (1989).
Following earlier work of Grimmett and Newman (1990) on the direct product
of a regular tree and $\Z$,
a general study of percolation on discrete groups
was initiated in Benjamini and Schramm (1996b).
One phenomenon that was
conjectured there to be general was a converse to the Burton and Keane
result, namely, that on any nonamenable group,
for some $p$, there are a.s.\ infinitely many
infinite clusters in $\omega_p$.
This led to the definition
$$
p_u(G)
:= \inf\Big\{p \st
   \P[\hbox{$\omega_p$ has exactly one infinite cluster}]=1
\Big\}
\,.
$$
Thus, $p_u(G) = p_c(G)$ when $G$ is an amenable Cayley graph.
H\"aggstr\"om and Peres (1997) show that on a Cayley
graph $G$, for {\it every\/} $p>p_u$,
there is exactly one infinite cluster a.s.\ in $p$-Bernoulli percolation.
It is known that $p_u < 1$ in many cases besides amenable groups,
e.g., finitely presented groups
with one end (Babson and Benjamini 1998) and Kazhdan groups (Lyons and
Schramm 1998).

The {\bf uniqueness phase} of Bernoulli percolation is the
range of $p$ where there is precisely one infinite cluster a.s.;
the {\bf nonuniqueness phase} is the range of
$p$ where there is more than one infinite cluster a.s.

The unique infinite percolation cluster of supercritical 
Bernoulli percolation on a graph, if there is such, can be viewed as a random 
perturbation of the graph. It is then natural to ask which properties of
the graph are inherited by such a percolation cluster.  

After presenting further 
definitions and reviewing some background in \ref s.back/,
we begin by studying in \ref s.geom/ purely geometric aspects of
percolation clusters, namely, how the isoperimetric constant $\isoe(G)$
(see \ref s.back/) behaves under percolation.
If $\omega$ is a random configuration of
Bernoulli percolation on a Cayley graph $G$,
then, of course, $\isoe(\omega) = 0$ a.s. 
However,

\procl t.isoeBern
If $G$ is a nonamenable Cayley graph and $\omega$ is a random configuration
of Bernoulli percolation
on $G$, then a.s.\ every infinite cluster of $\omega$ contains a subgraph
$\omega'$ with $\isoe(\omega') > 0$.
\endprocl

\noindent
(See \ref t.geom/.)
In fact, we show in \ref t.usubfr/ that one can further require $\omega'$ 
to be a tree. Note that $\omega$ is a random variable, respresenting the
configuration; in the sequel, we shall often say, however, that $\omega$ is
a ``percolation''.

The theorem raises the following question about percolation in
$\Z^d$.  Recall that the {\bf isoperimetric dimension} of a graph $G$
is the supremum of all $s$ such that 
$$
\inf\left\{ {|\bde V_1|^{s} \over |V_1|^{s-1}}\st V_1\subset\verts(G),\ 
0<|V_1|<\infty \right\}>0
\,.
$$

\procl q.zdisop
Let $\omega$ be supercritical Bernoulli percolation in $\Z^d$.
Is it true that for every $\epsilon>0$, a.s.\ $\omega$ contains a subgraph
with isoperimetric dimension at least $d-\epsilon$?  Does it contain a subgraph
with isoperimetric dimension $d$?
\endprocl

It is well known that if $G$ is
nonamenable, then the speed of simple random walk $\Seq{X(t)}$ on $G$
is positive, i.e., 
$$
\lim_{t\to\infty}{\distsub_G\big(X(0),X(t)\big)\over t}>0
\quad\hbox{a.s.}
\,,
$$
which results from the fact that the spectral radius is less than 1
(Kesten 1959); see \ref s.back/.
We prove the following extension
in \ref s.speed/ as a consequence of (a more precise version of)
\ref t.isoeBern/:

\procl t.speedBern
Let $G$ be a nonamenable Cayley graph and $\omega$ be Bernoulli percolation
on $G$. Let $\Seq{X(t)}$ be simple random walk on $\omega$. Given that the
cluster $K$ of $X(0)$ is infinite, we have a.s.\ that 
the speed of $X$ is positive.
\endprocl

\noindent
(See \ref t.speed/.)

We conjecture the following generalization:

\procl g.speedBern 
If $G$ is a Cayley graph on which simple random walk has positive speed,
then 
a.s., simple random walk on each infinite cluster of
$p$-Bernoulli percolation has positive speed.
\endprocl

There are (nontransitive) graphs
on which simple random walk has zero speed,
but for which Bernoulli percolation a.s.\ produces clusters
where simple random walk has positive speed.
An example is a binary tree with
a copy of $\Z$ attached to every vertex.

\procl g.psp %
If $G$ is a Cayley graph on which simple random walk has
zero speed, then a.s., simple random walk on every cluster of
Bernoulli percolation also has zero speed.
\endprocl

\noindent
Any possible counterexample
would have to be amenable and of exponential
growth by \ref l.SEBperc/ below.

The speed of random walk is related to other probabilistic behavior through
the following theorem due to the work of Avez (1974),
Derriennic (1980), Kaimanovich and Vershik (1983), and Varopoulos (1985).
Recall that a function $F : \verts(G) \to \R$ is
called {\bf harmonic} if $F(x) = \sum_{y \sim x} F(y)/\deg_G  x $ for all $x
\in \vertex$.

\procl t.SEB
\procname{Speed, Entropy, and Bounded Harmonic Functions}
The following conditions are equivalent for a given Cayley graph:
\beginitems
\itemrm{(i)}
the speed of simple random walk is zero;
\itemrm{(ii)}
the asymptotic entropy of simple random walk is zero;
\itemrm{(iii)}
there are no nonconstant bounded harmonic functions.
\enditems
\endprocl

Furthermore, Kaimanovich (1990) extended the equivalence of (ii) and (iii)
to many random walks in a random environment (RWRE)
that have a stationary measure; the extension to include (i)
is easy (\ref l.SEBperc/).
Since simple random walk restricted to percolation clusters has an
equivalent invariant measure (see \ref l.statmsr/), our
conjectures and results about the speed of random walk have some
alternative formulations in terms of entropy and bounded harmonic
functions.

Transience holds more generally, of course, than positive speed.
Transience of the infinite clusters of
Bernoulli percolation in $\Z^d$, $d > 2$, 
was established in Grimmett, Kesten and Zhang (1993) (see Benjamini,
Pemantle and Peres (1998) for a different proof). 
De Masi, Ferrari, Goldstein and Wick (1989)
proved an invariance principle for simple random walk on the
supercritical cluster in $\Z^2$.

\procl g.traBern If $G$ is a transient Cayley graph, then a.s.\ every
infinite cluster of Bernoulli
percolation on $G$ is transient.
\endprocl

The nonamenable case follows from \ref t.speedBern/; a slight extension is:

\procl t.gtraBern 
Let $G$ be a Cayley graph such that the ball of radius $n$ has cardinality
at least $\zeta^n$ for all $n$, where $\zeta>1$.
Then a.s., every infinite cluster of $p$-Bernoulli
percolation on $G$ is transient when $p > 1/\zeta$.
\endprocl

\noindent (See \ref t.gtra/.)

In \ref s.HD/, we study
the existence of
nonconstant harmonic Dirichlet functions
on the percolation clusters of Cayley graphs. 
A function $F : \verts(G) \to \R$ is
called {\bf Dirichlet} if $\sum_x \sum_{y \sim x} |F(x) -
F(y)|^2 < \infty$.
Recall that $\OHD$ denotes the class of graphs that do not
admit any nonconstant harmonic Dirichlet functions.
As we shall see below, the picture in the nonamenable situation is rather
involved and our understanding is far from complete.
In the amenable case, we have

\procl t.amenBern
If $G$ is an amenable Cayley graph
and $\omega$ is Bernoulli percolation,
then a.s.\ every cluster of $\omega$ is in $\OHD$.
\endprocl

\noindent (See \ref t.amen/.)

Medolla and Soardi (1995) proved that 
amenable transitive graphs are in $\OHD$ (see Remark~7.5 of \BLPSusf\
for a short proof via uniform spanning forests).
Soardi (1993) proved that $\OHD$ is invariant under quasi-isometries.
(See Lyons and Peres (1998) for a simple proof due to Schramm.)

Our proof of \ref t.amenBern/
uses the uniform 
spanning forest measures and their connection to harmonic Dirichlet
functions as presented in \BLPSusf. 
Such an ingredient can be motivated as follows:
In order to study the influence of 
geometric properties on potential-theoretic 
behavior, it is useful to have a geometric representation of the
potential-theoretic objects. The uniform spanning forest achieves this by
representing the analytic space of harmonic Dirichlet functions by a 
random geometric object.
The relevant definitions and properties are given in \ref s.HD/.

A positive answer to the following question would extend \ref t.amenBern/:

\procl q.ohd
Let $G$ be a Cayley graph, and suppose that $G\in\OHD$.
Let $\omega$ be Bernoulli percolation on $G$ in the uniqueness phase.
Does it follow that a.s.\ the infinite cluster of $\omega$ is in $\OHD$?
We do not know the answer even in the case that $G$ is the direct product
of a tree and $\Z$ or is a lattice in hyperbolic space $\H^d$ ($d \ge 3$).
\endprocl

In the nonuniqueness phase of Bernoulli percolation
on a (nonamenable) Cayley graph $G$,
the infinite clusters are not in $\OHD$ (\ref c.nuqHD/).

In the other direction,
for graphs admitting nonconstant harmonic Dirichlet functions, we believe:

\procl g.nhd
Let $G$ be a Cayley graph, $G\notin\OHD$.
Then a.s.\ all infinite clusters of $p$-Bernoulli percolation
are not in $\OHD$.
\endprocl

We can prove this for $p$ sufficiently large:

\procl t.notOHDBern
Let $G$ be a Cayley graph.
If $G \notin \OHD$, then there is some $p_0 < 1$ such that for every $p \ge
p_0$, almost surely no infinite
cluster of $p$-Bernoulli percolation
is in $\OHD$.
\endprocl

\noindent (See \ref t.notOHD/.)

The last section of our paper discusses questions
concerning the speed of simple random walk 
on a graph and a variant of the isoperimetric constant, called {\bf anchored
expansion}. The  anchored expansion constant
might be useful in studying the speed of 
simple random walk on infinite percolation clusters and other graphs as 
well.

Although in the above theorems, only Cayley graphs are
mentioned, we work in the greater generality of transitive graphs.
Similarly, we discuss percolation processes that are much more
general than Bernoulli percolation.  

\ssection{Acknowledgments}
We thank the organizers of the Cortona conference, 
Vadim A.\ Kaimanovich, Massimo A.\ Picardello,
Laurent Saloff-Coste, and Wolfgang Woess,
for a most enjoyable and productive meeting.
We thank Yuval Peres for comments on an earlier version
of this paper.

\bsection{Notation and Background}{s.back}

\ssection{Graph terminology, isoperimetric constant %
 and ends}
We use the letter $G$ to denote a graph and $\Gamma$ to denote
a closed subgroup of the automorphism group $\Aut(G)$ of $G$.
The vertices and edges of $G$ will be denoted
$\verts(G)$ and $\edges(G)$, respectively.
When there is an edge in $G$ joining vertices $u,v$, we write $u\sim v$.
The {\bf degree} $\deg v =\deg_G v $ of a vertex $v\in\verts(G)$ is the
number of edges incident with it.
A {\bf tree} is a connected graph with no cycles.
A {\bf forest} is a graph whose connected components are
trees.
The {\bf distance} between two vertices $v,u\in\verts(G)$
is denoted by $\dist(v,u)=\dist_G(v,u)$, and is the least number of edges
of a path in $G$ connecting $v$ and $u$. 

For a set of vertices $V_1\subset \verts(G)$,
let $\bde V_1$ denote the set of edges in $\edges(G)$ that
have one endpoint in $V_1$ and one endpoint in $\verts(G)\setminu V_1$.
The graph $G$ is {\bf amenable} if there is a sequence
of finite vertex subsets $V_1\subset V_2 \subset \cdots \subset
V_n\subset\cdots\subset\verts(G)$, such that $\bigcup_n V_n=\verts(G)$
and $|\bde V_n|/|V_n|\to 0$ as $n\to\infty$.
Here and in the sequel, $|A|$ denotes the cardinality of
a set $A$.
The (edge) {\bf isoperimetric constant} of a graph $G$,
also known as the {\bf Cheeger constant}, is defined by
$$
\isoe(G) :=\inf\bigl\{|\bde V_0|/|V_0|\st
\emptyset \ne V_0\subset \verts(G),\, |V_0|<\infty
\bigr\}
\,.
$$

An infinite set of vertices $V_0\subset\verts(G)$
is {\bf end-convergent} if for every finite $K\subset\verts(G)$,
there is a component of $G\setminu K$ that contains all
but finitely many vertices of $V_0$.
Two end-convergent sets $V_0,V_1$ are {\bf equivalent}
if $V_0\cup V_1$ is end-convergent.
An {\bf end} of $G$ is an equivalence class of end-convergent sets.

\ssection{Spectral radius, speed and entropy} 
Given $v,u\in\verts(G)$, let $p_t(v,u)$ be the
probability that simple random walk
starting at $v$ will be at $u$ at time $t$.
The {\bf spectral radius} $\rho(G)$ is defined by
$$
\rho(G) := \limsup_{t\to\infty} p_t(v,u)^{1/t}
$$
and does not depend on the choice of $v$ and $u$.
Dodziuk's (1984) discrete version of Cheeger's inequality 
states that if $G$ has bounded degrees and $\isoe(G)>0$,
then $\rho(G)<1$.

The {\bf speed} of a random walk $X$ starting at $o\in\verts(G)$ is
$$
\Lambda = \Lambda(X) := 
\lim_{t\to\infty}\dist\big(o,X(t)\big)/t
$$
when the limit exists.  The $\liminf$ speed is defined by
$$
\Lambda^- = \Lambda^-(X) := 
\liminf_{t\to\infty}\dist\big(o,X(t)\big)/t
\,.
$$
If $G$ has bounded degrees and $\rho(G)<1$,
then there are constants $\zeta>0$ and $\beta<1$
such that
$$
\PBig{\dist\bigl(o,X(t)\bigr)<\zeta t} \leq \beta^t
\,,
$$
because the probability that the random walk is inside
the ball of radius $\zeta t$ about $o$ is bounded
by the number of vertices in the ball times the probability that
$X$ is at the most likely vertex.
Consequently, in this situation, $\Lambda^->0$ a.s.

The {\bf entropy} of a probability measure $\mu$ on a finite or
countable set $A$
is defined to be
$$
H(\mu) := \sum_{x \in A} -\mu(x) \log \mu(x)
\,.
$$
Let $\mu_t$ denote the distribution of the location
$X(t)$ of a random walk $X$ at time $t$.
If $\lim_t H(\mu_t)/t$ exists, it is called the {\bf asymptotic entropy} of
the random walk. If $\Seq{X(t)}$ is simple random walk on a Cayley graph, then
the asymptotic entropy exists. %
In fact, the Subadditive Ergodic Theorem ensures the existence of 
$\lim_t -t^{-1} \log \mu_t\bigl(X(t)\bigr)$ a.s.\ and in $L^1$; see
Derriennic (1980).
The same is true for stationary RWRE, as observed by Kaimanovich (1990).
Similar reasoning applies to graphs with a transitive unimodular
automorphism group (Kaimanovich and Woess 1998).

For further information about random walks, spectral radius,
harmonic functions, etc., see Kaimanovich and Vershik (1983) and Woess (1994).

\ssection{Automorphism groups, unimodularity, and the Mass-Transport Principle}
Let $\Gamma\subset\Aut(G)$ be a subgroup of automorphisms of $G$ with the
topology of pointwise convergence.
We say that $\Gamma$ is (vertex) {\bf transitive} if for every
$v,u\in\verts(G)$,
there is a $\gamma\in\Gamma$ with $\gamma u=v$.
The graph $G$ is {\bf transitive} if $\Aut(G)$ is transitive.
Recall that every closed subgroup $\Gamma\subset \Aut(G)$ has a unique
(up to a constant scaling factor) Borel measure that, for every
$\gamma\in\Gamma$,
is invariant under left multiplication by $\gamma$; this 
measure is called (left) {\bf Haar measure}.
The group $\Gamma$ is {\bf unimodular} if Haar measure is also invariant
under right multiplication.

Most of our theorems concern percolation that is invariant
under a transitive unimodular closed subgroup of $\Aut(G)$.  
For example, when $G$ is the (right) Cayley graph of $\Gamma$ and
$\Gamma$ acts by left multiplication, then $\Gamma\subset\Aut(G)$
is (obviously) closed, unimodular, and transitive.
If $G$ is an amenable graph, then every closed transitive
subgroup of $\Aut(G)$ is unimodular (Soardi and Woess 1990). 

Several illustrations of the significance of unimodularity can
be found in \BLPSgip.
The most important one seems to be that when $\Gamma\subset\Aut(G)$
is unimodular, the Mass-Transport Principle 
takes the following simple form:

\procl t.mtp \procname{Mass-Transport Principle}
\assumptions.
Let $o\in\verts(G)$ be an
arbitrary basepoint.  Suppose that
$\phi:\verts(G)\times\verts(G)\to [0,\infty]$
is invariant under the diagonal action of $\Gamma$.
Then
$$
\sum_{v\in\verts(G)} \phi(o,v) = \sum_{v\in\verts(G)} \phi(v,o)
\,.
\label e.mtp
$$
\endprocl

See \BLPSgip\ for a discussion of this principle and for a proof.
In fact,  $\Gamma$ is unimodular iff \ref e.mtp/ holds for
every such $\phi$.  Hence, \ref e.mtp/ can be taken as
a definition of unimodularity.

\ssection{Percolation terminology}
A {\bf bond percolation} $\omega$ on $G$ is a random subset of $\edges(G)$.
For a more precise definition,
given a set $A$, let $2^A$ be
the collection of all subsets $\eta\subset A$, equipped with the
$\sigma$-field generated by the 
events $\{a\in\eta\}$, where $a\in A$.
A bond percolation $\omega$
on $G$ is then a random variable whose distribution is
a probability measure $\P$ on $2^{\edges(G)}$.
Similarly, a {\bf site percolation} is given by a probability
measure on $2^{\verts(G)}$, while a (mixed) {\bf percolation}
is given by a probability measure on $2^{\verts(G)\cup\edges(G)}$
that is supported on subgraphs of $G$.
If $\omega$ is a bond percolation, then $\hat\omega :=\verts(G)\cup\omega$
is the associated mixed percolation.
In this case, we shall often not distinguish between $\omega$ and
$\hat\omega$, and think of $\omega$ as a subgraph of $G$.
Similarly, if $\omega$ is a site percolation, there is an
associated mixed percolation
$\hat\omega:=\omega\cup\bigl(\edges(G)\cap(\omega\times\omega)\bigr)$,
and we shall often not bother to distinguish between $\omega$ and
$\hat\omega$.

Let $p\in[0,1]$.  Then the distribution of $p$-{\bf Bernoulli bond percolation}
$\omega$ on $G$ is the product measure on $2^{\edges(G)}$
that satisfies $\P[e\in\omega]=p$ for all $e\in\edges(G)$.
Similarly, one defines $p$-{\bf Bernoulli site percolation} on
$2^{\verts(G)}$.

If $v\in\verts(G)$ and $\omega$ is a percolation on $G$,
the {\bf component} (or {\bf cluster}) $K(v)$ of $v$ in $\omega$
is the set of vertices in $\verts(G)$ that can be connected
to $v$ by paths contained in $\omega$.

Suppose that $\Gamma$ is an automorphism group of a graph $G$.
A percolation on $G$ is {\bf $\Gamma$-invariant}
if its distribution $\P$ is invariant under each automorphism in $\Gamma$.

\ssection{Insertion tolerance and component indistinguishability}
Given a set $Z$, an element $z\in Z$, and a subset $\omega\in 2^Z$, 
let $\ins_z\omega:=\omega\cup\{z\}$.
A probability measure $\P$ on $2^Z$
is {\bf insertion tolerant} if there is a constant $\delta>0$ such that
$\P[\ins_z\ev A]\geq \delta\P[\ev A]$
holds for each measurable $\ev A \subset 2^Z$ and
each $z\in Z$.
Loosely, this means that inserting any $z$ into $\omega$ does
not decrease its probability by more than a constant factor.
For example, $p$-Bernoulli bond percolation is insertion tolerant
when $p>0$.

Let $G$ be graph and $\Gamma$ a closed transitive subgroup
of $\Aut(G)$.  Let $\omega$ be a $\Gamma$-invariant bond percolation.
We say that $\omega$ has {\bf indistinguishable components}
if for every measurable
$\ev A\subset 2^{\verts(G)}\times 2^{\edges(G)}$
that is invariant under the diagonal action of $\Gamma$,
almost surely, for all infinite components $C$ of $\omega$,
we have $(C,\omega)\in\ev A$, or for all
infinite components $C$, we have $(C,\omega)\notin\ev A$.
(That is, whether $(C, \omega) \in \ev A$ does not depend on $C$, but may
depend on $\omega$.)

The following is from Lyons and Schramm (1998):

\procl t.cerg \procname{Component Indistinguishability} 
\assumptions.
Every $\Gamma$-invariant, insertion tolerant, bond percolation process on $G$
has indistinguishable components.
\endprocl

For example, this shows that in $p$-Bernoulli 
bond percolation on a Cayley graph of a nonamenable group,
almost surely, either all infinite clusters are
transient, or all clusters are recurrent.
In fact, as indicated by \ref t.speedBern/, a.s.\ all infinite clusters
are transient.

Similar statements hold for site and mixed percolations.

\bsection{Geometry of Perturbations of Nonamenable Graphs}{s.geom}

Let $G$ be an infinite graph and $K$ a finite subgraph
of $G$.  Set
$$
\alpha_K:= {1 \over |\verts(K)|} \sum_{v\in\verts(K)} \deg_K v 
\,,
$$
where $\deg_K v$ is the degree of $v$ in $K$.
Define
$$
\alpha(G) := \sup\bigl\{ \alpha_K \st K\subset G\hbox{ is finite}\bigr\}
\,.
$$
Note that  when all vertices in $G$ have degree $d$, the
isoperimetric constant of $G$ satisfies
$$
\isoe(G) = d - \alpha(G)
\,.
$$

Let $T$ be a regular tree and $o\in\verts(T)$ be
some basepoint.
H\"aggstr\"om (1997) has shown that when $\omega$
is an automorphism-invariant percolation
on  $T$ and $\E[\deg_\omega o]\geq\alpha(T)$, there are
infinite clusters in $\omega$ with positive probability.
In \BLPSgip, it was shown that the same
result applies to transitive graphs with a unimodular
automorphism group.  In \ref t.subch/, we extend this result
and show that with the same assumptions
but with a strict inequality $\E[\deg_\omega o]>\alpha(G)$, with positive
probability, there is a subgraph in $\omega$ with
$\isoe>0$.
This will be used to prove

\procl t.usubgr \procname{Uniqueness Gives a Subgraph with $\isoe>0$}
\assumptions, and suppose that $\isoe(G)>0$.
Let $\omega$ be a $\Gamma$-invariant percolation
in $G$ that has a.s.\ exactly one infinite component.
Then (on a larger probability space) there is
a percolation $\omega'\subset\omega$ 
such that $\omega'\neq\emptyset$ and $\isoe(\omega')>0$  a.s.
Moreover, the distribution of the
pair $(\omega',\omega)$ is $\Gamma$-invariant.
\endprocl

In the following, if $K\subset G$ is a subgraph and
$v\in\verts(G)$ is not in $K$, then we set $\deg_K v:=0$.

As we indicated,
the proof of \ref t.usubgr/ is based on the following
more quantitative result.

\procl t.subch \procname{High Marginals Give a Subgraph with $\isoe>0$}
\assumptions. 
Let $\omega$ be a $\Gamma$-invariant
(nonempty) percolation in $G$.  Let $h>0$ and suppose that
$$
\E[\deg_{\omega}o \mid o\in\omega] > \bigl( \alpha(G)+ 2 h\bigr)
\,.
\label e.chbd
$$
Then there is (on a larger probability space)
a percolation $\omega'\subset\omega$ such
that $\omega'\neq\emptyset$ and $\isoe(\omega')\geq h$ 
with positive probability.
Moreover, the distribution of the pair $(\omega',\omega)$ is
$\Gamma$-invariant.
\endprocl

\proof  Given any subgraph $\omega$ of $G$,
we define percolations $\omega_n$ on $\omega$ inductively 
as follows.  Set $\omega_0:=\omega$.
Suppose that $\omega_n$ has been defined.
Let $\beta_n$ be a $(1/2)$-Bernoulli site percolation on
$G$, independent of $\omega_0,\dots,\omega_n$.
Let $\gamma_n$ be the union of the finite components $K$ of
$\beta_n \cap \omega_n$ that satisfy %
$$
{|\partial_{\edges(\omega_n)} K| \over |K|} < h
\,,
$$
where $\partial_{\edges(\omega_n)} K$ denotes the set of edges of
$\omega_n$ connecting $K$ to its complement.
Now set $\omega_{n+1}:=\omega_n\setminu \gamma_n$.
Finally, define 
$$
\omega':=\bigcap_{n=0}^\infty \omega_n
\,.
$$
For future use, write $\Xi(h, \omega) := \omega'$.

Most of the proof will be devoted to showing that
$\omega'\neq\emptyset$ with positive probability,
but first we verify that $\isoe(\omega')\geq h$.
Indeed, let $W$ be a finite nonempty set of vertices
in $G$, let $F$ be the set of all edges
of $G$ incident with $W$, and let $F_0\subset F$.
Suppose that $|F_0|/|W|<h$.
To verify that $\isoe(\omega')\geq h$ a.s.,
it is enough to show that the probability that
$W\subset \omega'$ and $\omega'\cap F=F_0$ is zero.  
If $\omega_n\cap F = F_0$ for some $n$, then a.s.\
there is some $m>n$ such that $W$ is a component
of $\beta_m$. 
Now either $W \not\subset \omega_m$, in which case $W \not\subset \omega'$,
or $W \subset \omega_m$, in which case
$W$ is not contained in $\omega_{m+1}$,
hence not in $\omega'$.  On the other hand,
if $\omega_n\cap F\neq F_0$ for every $n$,
then also $\omega'\cap F\neq F_0$.  Consequently
$\isoe(\omega')\geq h$ a.s.

Now set 
\begineqalno
D_n & := \E \deg_{\omega_n} o\,,
\qquad
D_\infty  := \E \deg_{\omega'} o\,,
\cr
\theta_n &:= \P[o\in\omega_n]\,,
\qquad
\theta_\infty := \P[o\in\omega']
\,.
\endeqalno
Our goal is to prove the inequality
$$
D_{n+1}\geq D_n - (\theta_n-\theta_{n+1})\bigl(\alpha(G)+2 h\bigr)
\,.
\label e.dec
$$
This will be achieved through use of the Mass-Transport Principle.
Observe that
$$
\theta_n-\theta_{n+1} = \P[o\in\gamma_n]
\,.
$$
Fix $n$ and
define the random function $m:\verts(G)\times\verts(G)\to[0,\infty)$
as follows.
For every vertex $v\in\verts(G)$, let $K(v)$ be
the component of $v$ in $\gamma_n$, which we take to be $\emptyset$
if $v\notin\gamma_n$.
Let $v,u\in\verts$.
If $u\notin\gamma_n$, set $m(v,u):=0$. 
If $v\in K(u)$, let $m(v,u):=\deg_{\omega_n}v/|K(u)|$.
Otherwise, let $m(v,u)$ be $|K(u)|^{-1}$ times the number of edges in
$\omega_n$ that connect $v$ to a vertex in $K(u)$.
Note that $v$ and $u$ need not be adjacent in order that $m(v, u) \ne 0$.
Clearly, $\E m(v,u)$ is invariant under the diagonal
action of $\Gamma$ on $\verts(G)\times\verts(G)$.
Consequently, the \MTP/ implies that
$$
\sum_{v\in\verts(G)}\E m(o,v) =\sum_{v\in\verts(G)}\E m(v,o)
\,.
$$
A straightforward calculation shows that
$$
\sum_{v\in\verts(G)} m(o,v) = \deg_{\omega_n}o-\deg_{\omega_{n+1}}o
\,,
$$
while, if $o\in\gamma_n$, we have that
$\sum_{v\in\verts(G)} m(v,o)$ is equal to twice
the number of edges of $\omega_n$ incident with $K(o)$, divided
by $|K(o)|$.  The
number of edges of $G$ with both endpoints in $K(o)$
is at most $\alpha(G)|K(o)|/2$, and,
by construction, $|\partial_{\omega_n}K(o)|< h|K(o)|$.
Hence
$$
\sum_{v\in\verts(G)} m(v,o) < \alpha(G) + 2 h
\label e.alh
$$
when $o\in\gamma_n$ and $ \sum_{v\in\verts(G)} m(v,o)=0$ otherwise.
Therefore, 
\begineqalno
D_n - D_{n+1}
&=
\E[\deg_{\omega_n}o-\deg_{\omega_{n+1}}o]
=
\sum_{v\in\verts(G)}\E m(v,o)
\cr&\le
(\alpha(G) + 2 h) \P[o\in\gamma_n]
=
(\alpha(G) + 2 h) (\theta_n-\theta_{n+1})
\,,
\cr
\endeqalno
which is the same as \ref e.dec/.

Induction and \ref e.dec/ give
$$
D_n \geq D_0 - \theta_0\bigl(\alpha(G)+2 h\bigr) +
\theta_n\bigl(\alpha(G)+2 h\bigr)
\,;
$$
taking a limit as $n\to\infty$ yields the inequality
$$
D_\infty \geq D_0 - \theta_0\bigl(\alpha(G)+2 h\bigr) +
\theta_\infty\bigl(\alpha(G)+2 h\bigr)
\,.
\label e.fin
$$
This gives $D_\infty>0$, because
 \ref e.chbd/ is
equivalent to $D_0 - \theta_0\bigl(\alpha(G)+2 h\bigr)>0$.
Consequently, $\omega'\neq\emptyset$ with positive
probability.
\Qed

\procl r.est \procname{The Density of $\omega'$}
The following lower bound for $\theta_\infty$ is a
consequence of
\ref e.fin/ and the inequality $\theta_\infty\deg_G o \geq D_\infty$:
\begineqalno
\P[o\in\omega'] = \theta_\infty
& \geq
{D_0-\bigl(\alpha(G)+2 h\bigr)\theta_0
\over
\deg_G o - \bigl(\alpha(G)+2 h\bigr)}
\cr &
=
\P[o\in\omega]
\left(
1-
{\deg_G o - \E[\deg_{\omega}o \mid o\in\omega] \over \isoe(G) - 2 h}
\right)
\,.
\label e.denlb
\endeqalno
\endprocl

\procl r.weakiq \procname{A Weak Inequality Suffices}
In fact, in place of \ref e.chbd/, it is enough
to assume the weak inequality
$ \E[\deg_{\omega}o \mid o\in\omega] \geq \bigl( \alpha(G)+ 2 h\bigr)$.
The reason is that
the inequality \ref e.alh/ is strict when
$o\in\gamma_n$, which implies that \ref e.fin/ is strict when
$D_\infty\neq D_0$.
\endprocl

\procl t.subtr \procname{Threshold for a Forest}
If $\omega$ is a forest a.s., then \ref t.subch/
is true when $\alpha(G)$ is replaced by $2$.
\endprocl

\proof
In a finite tree $K\subset G$, we have
$\alpha_K<2$.  Hence the proof of \ref t.subch/ applies
with $2$ replacing $\alpha(G)$  everywhere.
\Qed

\proofof t.usubgr
Fix a basepoint $o\in\verts(G)$.
Let $\omega_*$ be the infinite component of $\omega$.
Conditioned on $\omega$,
for every vertex $v\in\verts(G)$, let $\phi(v)$
be chosen uniformly among the vertices of
$\omega_*$ closest to $v$, with all $\phi(v)$ independent given $\omega$,
and for edges $e=[v,u]\in\edges(G)$, let $\phi(e)$
be chosen uniformly among shortest paths in $\omega_*$
joining $\phi(v)$ and $\phi(u)$, with 
all $\phi(e)$ independent given all $\phi(v)$ and $\omega$.
For integers $j$, let $\eta_j$ be the set of edges
$e\in\edges(G)$ such that $\phi(e)$ is contained
within a ball of radius $j$ about one of the
endpoints of $e$.
Then $\eta_1\subset\eta_2\subset\cdots$ are 
$\Gamma$-invariant bond percolations
on $G$ with $\bigcup_j\eta_j =\edges(G)$.
Consequently, $\E \deg_{\eta_j}o\to\deg_G o$
as $j\to\infty$.
For each $j$, choose independently a random sample of $\Xi(h, \eta_j)$ and
denote it $\xi_j$.
By \ref e.denlb/, we have that $\P[\xi_j \ne \emptyset] \to 1$.
Let $J := \inf\{j \st \xi_j \ne \emptyset\}$. Then $J <
\infty$ a.s.
Set $\omega' := \phi\left(\xi_{J}\right)$.
Since $\isoe(\xi_J) \ge h$ a.s., we have also
$\isoe\bigl(\phi(\xi_J)\bigr)>0$ a.s.
\Qed

Suppose that $G$ is transitive, $\isoe(G)>0$, and
$\omega$ is, say, Bernoulli percolation on $G$ that has
a.s.\ more than one infinite component. 
Then \ref t.usubgr/ does not apply to $\omega$.
However, as observed by Burton and Keane (1989), insertion tolerance shows
that there
are a.s.\ components of $\omega$ with at least three ends.
Hence the next theorem does apply.

\procl t.nuch
\procname{A Forest with $\isoe>0$ Inside Many-Ended Percolation}
\assumptions, and let $\omega$ be a $\Gamma$-invariant
percolation on $G$.  Suppose that a.s.,
there are components of $\omega$ with at least three ends.
Then there is (on a larger probability space) a
random forest $\fo\subset\omega$ with $\isoe(\fo)>0$,
$\fo\neq\emptyset$ a.s.,
and the distribution of the pair $(\fo,\omega)$
is $\Gamma$-invariant.
\endprocl

We shall need the following two lemmas from \BLPSgip.
Recall that $K(x)$ denotes the component of $x$ in $\omega$.

\procl l.ends \procname{Ends, $p_c$ and Degrees}
\assumptions.
Let $\omega$ be a $\Gamma$-invariant percolation on $G$ that
has infinite components with positive probability.
If
\beginitems
\itemrm{(i)} some component of $\omega$ has at least three ends with
positive probability,
\enditems
then
\beginitems
\itemrm{(ii)} some component of $\omega$ has $p_c < 1$ with positive
probability and
\itemrm{(iii)} for every vertex $x$, $\Eleft{\deg_\omega x  \bigm|
|K(x)| = \infty} > 2$.
\enditems
If $\omega$ is a forest a.s., then the three conditions are equivalent.
\endprocl

\procl l.trim \procname{Trimming to a Forest}
\assumptions.
Let $\omega$
be a $\Gamma$-invariant percolation on $G$ such that a.s.\ there is a
component of $\omega$ with at least three ends.
Then (on a larger probability
space) there is a
random forest $\fo\subset\omega$ such that
the distribution of the pair $(\fo,\omega)$ is
$\Gamma$-invariant and a.s.\ 
whenever a component $K$ of $\omega$ has at least three ends, there
is a component of $K \cap \fo$ that has infinitely many ends.
\endprocl

\proofof t.nuch
By \ref l.trim/, there is a random forest $\fo'\subset\omega$
with  some components
having infinitely many ends a.s.\
and the distribution of $(\fo',\omega)$ is $\Gamma$-invariant.
Let $\fo''$ be the union of the infinite components of $\fo'$. 
By \ref l.ends/,
$\E[\deg_{\fo''}o \mid o\in\fo'']>2$.
Given $\omega$, for each $j$, let $\xi_j$ be an independent sample of $\Xi(1/j,
\fo'')$.
Put $\fo := \xi_J$,
where $J := \inf\{j \st \xi_j \ne \emptyset\} <
\infty$ a.s.: Clearly, $J < \infty$ with positive probability. If the set
$\ev A$ of $\omega$ where
$J = \infty$ had positive probability, then we would obtain a
contradiction to what has just been proved by noting that $\ev A$ is
$\Gamma$-invariant and by conditioning on $\ev A$.
\Qed

We can now deduce the following extension of \ref t.isoeBern/:

\procl t.geom
\assumptions, and suppose that $\isoe(G)>0$.
Let $\omega$ be a $\Gamma$-invariant percolation on $G$ that has infinite
clusters a.s.  
Then in each of the following cases
(on a larger probability space) there is
a percolation $\omega'\subset\omega$ 
such that $\omega'\neq\emptyset$, $\isoe(\omega')>0$  a.s., and
the distribution of the pair $(\omega',\omega)$ is $\Gamma$-invariant:
\beginitems
\itemrm{(i)} $\omega$ is Bernoulli percolation;
\itemrm{(ii)} $\omega$ has a unique infinite cluster a.s.;
\itemrm{(iii)} $\omega$ has a cluster with at least three ends a.s.;
\itemrm{(iv)} $\E[\deg_{\omega}o \mid o\in\omega] > \alpha(G)$ and
  $\omega$ is ergodic.
\enditems
\endprocl

\proof
In Bernoulli percolation, if there is more than one
infinite cluster, then there is a cluster with at least three ends
by insertion tolerance and ergodicity.  Consequently,
(i) follows from (ii) and (iii).  Parts (ii)--(iv)
follow from Theorems 
\briefref t.usubgr/, \briefref t.nuch/, and \briefref t.subch/.
\Qed

Although it will not be needed in the sequel,
we note that \ref t.usubgr/ can be strengthened as follows.

\procl t.usubfr \procname{A Forest in the Uniqueness Regime with $\isoe>0$}
\assumptions, and suppose that $\isoe(G)>0$.
Let $\omega$ be a $\Gamma$-invariant percolation
on $G$ that has a.s.\ exactly one infinite component.
Then (on a larger probability space) there is
a random forest $\fo\subset\omega$ with 
$\fo\neq\emptyset$ and
$\isoe(\fo)>0$ a.s., and the distribution of
the pair $(\fo,\omega)$ is $\Gamma$-invariant.
\endprocl

\proof Let $\omega'\subset\omega$ be as in \ref t.usubgr/.
Since $\isoe(\omega')>0$,
Theorem 13.7 from \BLPSusf\ constructs a percolation $\omega''\subset\omega'$
such that a.s.\ 
$\omega''$ has all components  with infinitely many ends and the
distribution of 
$(\omega'',\omega)$ is $\Gamma$-invariant. 
By \ref l.trim/, there is a forest
$\fo'\subset\omega''$ such that some tree in
$\fo'$ has infinitely many ends a.s.\
and the distribution of $(\fo',\omega)$ is $\Gamma$-invariant.
Let $\fo''$ be the union of the
infinite components of $\fo'$.
By \ref l.ends/, $\E[\deg_{\fo''}o \mid o\in\fo'']>2$.
The proof is completed as for \ref t.nuch/.
\Qed

Applying \ref t.usubfr/ to the case where $\omega=G$ a.s.,
we obtain an invariant random forest
in $G$ with $\isoe(\fo)>0$.
This is related to the result of
Benjamini and Schramm (1997) which says that every
bounded-degree graph with $\isoe>0$ contains
a tree $T$ with $\isoe(T)>0$.  
In fact, the latter result can be used to extend our
theory to the non-transitive setting as follows.

\procl c.ntrans \procname{The Non-Transitive Case}
Let $G$ be a graph of bounded degree with $\isoe(G)>0$.  Then there
is some $p_0<1$ such that $p$-Bernoulli bond
percolation on $G$ has a subgraph with $\isoe>0$ a.s.\ whenever $p>p_0$. 
\endprocl

\proof
By the result of Benjamini and Schramm (1997) mentioned above,
there is a tree $T\subset G$ with $\isoe(T)>0$.
Let $T_3$ be the $3$-regular tree.
There is a map $\phi$ that takes $\verts(T_3)$ into $\verts(T)$,
takes every edge $e=[v,u]\in\edges(T_3)$ to a path of
bounded length $\phi(e)$ in $T$ joining $\phi(v)$ to $\phi(u)$,
and when $e,e'\in\edges(T_3)$ are distinct, the corresponding paths
$\phi(e),\phi(e')$ are edge-disjoint.
Consequently, $p$-Bernoulli percolation on $T$ can be pulled back
via $\phi$ to a bond percolation $\omega$ on $T_3$ in which
the events $\{e\in\omega\}$ ($e\in\edges(T_3)$) are mutually independent.
Moreover, $\P[e\in\omega]\geq 1-k(1-p)$, where $k$
is the maximum length of a path $\phi(e')$, $e'\in\edges(T_3)$.
Consequently, $\omega$ dominates $\bigl(1-k(1-p)\bigr)$-Bernoulli bond
percolation $\omega'$ on $T_3$. 
By \ref t.subtr/, when $3\bigl(1-k(1-p)\bigr)>2$, there is 
with positive probability, and therefore a.s., a subgraph
$\omega''\subset\omega'$ with $\isoe(\omega'')>0$.
Now $\phi(\omega'')$ is the required subgraph of $\omega$.
\Qed

\bsection{Speed and Transience}{s.speed}

In this section, we prove that in many cases,
simple random walk on the infinite
components of invariant percolation on a nonamenable transitive
graph $G$ has positive speed.

Let $\omega$ be a percolation on $G$.
It will be useful to consider {\bf delayed simple random walk} $Z =
Z^\omega$ on $\omega$, defined as follows.
Set $Z(0):=o$, where $o\in\verts(G)$ is some fixed basepoint.
If $n\geq 0$, conditioned on $\Seq{Z(0),\dots,Z(n)}$ and $\omega$,
let $Z'(n+1)$ be chosen from $Z(n)$ and its
neighbors in $\edges(G)$ with equal probability.
Set $Z(n+1) := Z'(n+1)$ if the edge $[Z(n), Z'(n+1)]$ belongs to $\omega$;
otherwise, let $Z(n+1) := Z(n)$.

For each $x \in \verts(G)$ and $n\in\Z$,
choose $\phi_n(x)$ uniformly with respect to Haar
measure among all $\gamma\in\Gamma$
satisfying $\gamma x=o$. These choices are to be independent of each other,
of $\omega$, and of $Z$. For any sequence $\Seq{z(0),z(1),\dots}$,
let $\SS z$ be the shifted sequence
$\Seq{\SS z(0),\SS z(1),\dots}$ defined by $\SS z(n):=z(n+1)$.

The following is shown in Lyons and Peres (1998):

\procl l.statmsr
\procname{Stationarity of Random Walk}
\assumptions.
Let $\omega$ be a $\Gamma$-invariant percolation on $G$ with law $\mu$.
Set 
$$
\hat Z(n) := \bigl(\phi_n(Z(n)) \SS^n Z,\phi_n(Z(n)) \omega\bigr)
\,.
$$
The sequence $\hat Z$ is stationary
with respect to $\SS$
(in the big probability space where $\omega$ is also random).
\par
Let $\mu'$ be the measure on subgraphs $\omega\subset G$
whose Radon-Nikodym derivative with respect to $\mu$ is
$\deg_\omega o/\E_\mu[\deg_\omega o]$.
Let $X$ be simple random walk on $\omega$ starting at $o$,
and let $\hat X$ be defined analogously to $\hat Z$.
Then $\hat X$ is $\SS$-stationary when $\omega$ is chosen with the law
$\mu'$.
\endprocl

\procl l.nonrandom   \procname{Speed Exists and is Not Random}
\assumptions.
Let $\omega$ be a $\Gamma$-invariant percolation on $G$.
Then the speed $\Lambda$ of delayed simple random walk on $\omega$
exists and is an $\omega$-measurable random variable (possibly zero).
\par
If $\omega$ has indistinguishable components and is ergodic, then,
conditioned on $|K(o)|=\infty$, $\Lambda$ is equal a.s.\ to a constant.
\par
The same statements hold for simple random walk in place of delayed
simple random walk.
\endprocl

\proof
Let $f_n(\hat Z):=\distsub_G(o,Z(n))$.
Then 
$$
f_{n+m}(\hat Z) \leq f_n(\hat Z) + f_m(\SS^n\hat Z)
$$
by the triangle inequality.
Consequently, the Subadditive Ergodic Theorem shows that the speed
$$
\Lambda = \Lambda(\hat Z) = \lim_{n\to\infty} f_n(\hat Z)/n
=\lim_{n\to\infty}\distsub_G\big(o,Z(n)\big)/n
$$
exists a.s.

To show that
the speed $\Lambda(\hat Z)$ depends only on $\omega$ and not on the path of
the random walk a.s., 
define $F(\hat Z)$ to be the variance, conditioned on $\omega$,
of the speed of an independent random
walk starting from $Z(0)$. By L\'evy's 0-1 Law, $F(\SS^n \hat Z)$ converges
to zero a.s. But by stationarity, the distribution of $F(\SS^n \hat Z)$ is
the same for all $n$. Hence, it is 0.

The statement concerning indistinguishable components is a consequence of
the definition.

The same proof applies to simple random walk since the measures $\mu$ and
$\mu'$ of \ref l.statmsr/ are mutually absolutely continuous.
\Qed

Our main tool to convert the geometric information of \ref s.geom/
to probabilistic information is the following:

\procl t.sp \procname{Speed When There is a Subgraph with $\isoe>0$}
\assumptions. 
Let $\omega'\subset\omega$  be percolations on
$G$ such that the distribution of the pair $(\omega',\omega)$
is $\Gamma$-invariant.  Suppose that
$\omega'\neq\emptyset$ and $\isoe(\omega')>0$ a.s.
Then simple random
walk on $\omega$ has positive speed a.s.
\endprocl

\proof
Let $V^*$ be the vertices of $\omega'$ that are
in the $\omega$-component of $o$,
and
let $Z$ be delayed simple random walk on $\omega$ 
starting at $o$.
Note that given $\omega$,
$Z$ is reversible with uniform stationary distribution.
Given $\omega$ and $\omega'$ with $o\in \omega'$,
there is an induced walk $Z^*$ on $V^*$ defined as
follows.  Set $t_0:=0$. 
Since the transformed $\hat Z$ is stationary (\ref l.statmsr/),
the Poincar\'e recurrence theorem (see, e.g., Petersen (1983), p.~34)
shows that conditioned on $o\in \omega'$,
there is a.s.\ some first time $t_1>0$ such that
$Z(t_1)\in \omega'$.
(Strictly speaking, \ref l.statmsr/ does not apply
when there is an extra ``scenery" $\omega'$.
But the lemma extends easily to this situation.)
Inductively, for $k>0$, let $t_k$ be the
first time $t>t_{k-1}$ such that $Z(t)\in \omega'$.
Define $Z^*(k):=Z(t_k)$.  Then given $\omega$ and $\omega'$ with $o\in
\omega'$,
$Z^*$ is just the Markov chain $Z$ induced on the
states $V^*$.
In particular, it is reversible with the same stationary
distribution on $V^*$, i.e., uniform.

We claim that,
given $\omega$ and $\omega'$ with $o \in \omega'$,
the spectral radius $\rho(Z^*)$ is less than $1$ a.s.
Given two vertices $u^*,v^*\in V^*$, 
let $p^*(u^*,v^*)$ denote the transition probability
of the Markov chain $Z^*$.
Let $G^*$ be the graph whose vertices
are $V^*$ and whose edges $[u^*,v^*]$ are those pairs with
$p^*(u^*,v^*)>0$. 
Note that there is some positive lower bound $c>0$ for $p^*(u^*,v^*)$
whenever $[u^*,v^*]\in\omega'$.
Consequently,
$$
\inf \left\{
{1
\over
|K^*|}
\sum_{e^* \in \partial_{\edges} K^*} p^*(e^*)
\st
K^* \subset \verts(G^*) \hbox{ is finite}
\right\}
\geq c \isoe(\omega')
> 0
\,.
$$
Since the stationary distribution is uniform,
this implies that $\rho(Z^*)<1$ a.s.~(see, e.g., Kaimanovich (1992)),
as claimed.  

Fix some $\rho_0<1$ such that $\rho(Z^*)<\rho_0$ with positive
probability, and let $\ev A$ be the event that
$o\in\omega'$ and $\rho(Z^*)<\rho_0$.
Then, for some $\zeta<1$ and all $v\in\verts(G)$,
$$
\Pbig{Z(t_k)=v\,\bigm|\omega, \omega'}\leq \zeta^k \quad\hbox{on }\ev A
\,,
$$
which gives
$$
\Pbig{Z(t_k)=v\,\bigm|\ev A}\leq \zeta^k
\,.
\label e.vhit
$$
Since the number of vertices $v\in\verts(G)$ with
$\distsub_G(o,v)<r$ is bounded by $(\deg_G o)^{r+1}$, by summing
\ref e.vhit/ over all such vertices, we get
$$
\PBig{\distsub_G\big(o,Z(t_k)\big) <r\,\bigm| \ev A}
\leq \zeta^k (\deg_G o)^{r+1}
.
$$
Let $\beta$ be such that $(\deg_G o)^\beta = 1/\zeta$, and choose
$r:= (\beta/2)k-1$.
Then
$$
\PBig{\distsub_G\big(o,Z(t_k)\big) <\beta k/ 2 - 1\,\bigm| \ev A}
\leq \zeta^{k/2}
.
\label e.tkspeed
$$
By the Borel-Cantelli lemma, it follows that
$$
\liminf_{k \to\infty} \distsub_G\big(o,Z(t_k)\big)/k \ge \beta/2
$$
a.s.\ on $\ev A$.
Also, the ergodic theorem ensures that $\lim t_k/k < \infty$ a.s., whence
$$
\liminf_{k \to\infty} \distsub_G\big(o,Z(t_k)\big)/t_k > 0
$$
a.s.
This shows that the speed of $Z$ is positive a.s.\ by \ref l.nonrandom/.
By the obvious coupling of delayed random walk
and simple random walk, it follows that
also the speed of simple random walk is positive a.s.
\Qed

\procl t.speed \procname{Speed}
\assumptions, and suppose that $\isoe(G)>0$.
Let $\omega$ be a $\Gamma$-invariant percolation on $G$.
Then simple random walk on some
infinite cluster of $\omega$ has positive speed with positive
probability in each of the following cases:
\beginitems
\itemrm{(i)} $\omega$ is Bernoulli percolation that has
infinite components a.s.;
\itemrm{(ii)} $\omega$ has a unique infinite cluster a.s.;
\itemrm{(iii)} $\omega$ has a cluster with at least three ends with
positive probability;
\itemrm{(iv)} $\E[\deg_{\omega}o \mid o\in\omega] > \alpha(G)$.
\enditems
\endprocl

\proof
This follows from Theorems \briefref t.geom/ and \briefref t.sp/.
\Qed

In case $G$ is a tree, (iii) and (iv) of this theorem 
were established by H\"aggstr\"om (1997).

In case the percolation is ergodic and has indistinguishable components,
like Bernoulli percolation, we have the stronger conclusion that simple
random walk has positive speed on every infinite component a.s.

\procl r.sharpness
It does not suffice in 
\ref t.speed/ to drop in (i) the assumption that $\omega$ is Bernoulli.
For example, if $\omega$ is the wired uniform spanning forest
(see \ref s.HD/),
then every component is a tree with one end (\BLPSusf), whence is recurrent.
\endprocl

In order to derive additional consequences of \ref t.speed/, we now extend
\ref t.SEB/:

\procl l.SEBperc
\assumptions.
Let $\omega$ be a $\Gamma$-invariant bond percolation on $G$.
The following are equivalent:
\beginitems
\itemrm{(i)}
the speed of simple random walk $X(t)$ on $K(o)$ is zero a.s.\
in the $G$-metric:
$$
\lim_{t \to\infty} \dist_G\bigl(X(t)\bigr)/t = 0\,;
$$
\itemrm{(ii)}
the speed of simple random walk $X(t)$ on $K(o)$ is zero a.s.\
in the $\omega$-metric:
$$
\lim_{t \to\infty} \dist_\omega\bigl(X(t)\bigr)/t = 0
\,;
$$
\itemrm{(iii)}
the asymptotic entropy of simple random walk on $K(o)$ is zero a.s.;
\itemrm{(iv)}
there are no nonconstant bounded harmonic functions on $K(o)$ a.s.
\enditems
\endprocl

\proof Because of \ref l.statmsr/, the equivalence of
(iii) and (iv) follows from Kaimanovich and Woess (1998).
Clearly, (ii) implies (i). We show that (i) implies (iii) implies
(ii).

Assume (i). Fix $\omega$ such that the speed on $K(o)$ is zero.
Let $\mu^\omega_t$ denote the law of $X(t)$ on $K(o)$.
Let $B_r$ denote the ball of radius $\lfloor{r}\rfloor$ in $G$ centered at $o$.
Given $\epsilon > 0$, choose $t_0$ large enough that
for all $t \ge t_0$, we have
$\mu^\omega_t(B_{t\epsilon}) \geq 1-\epsilon$.
Let $D := \deg_G o$.
Then for $t \ge t_0$, concavity of $\log$ gives the inequality
$$
\sum_{x \in B_{t\epsilon}}
   -\mu^\omega_t(x) \log \mu^\omega_t(x)
\le
\mu^\omega_t(B_{t\epsilon}) \log \bigl(|B_{t\epsilon}| /
  \mu^\omega_t(B_{t\epsilon})\bigr)
\le
\log \bigl(D^{t\epsilon}/(1-\epsilon)\bigr)
\,.
$$
Similarly,
$$
\sum_{x \notin B_{t\epsilon}}
   -\mu^\omega_t(x) \log \mu^\omega_t(x)
=
\sum_{x \in B_t\setminu B_{t\epsilon}}
   -\mu^\omega_t(x) \log \mu^\omega_t(x)
\le
\epsilon \log (D^{t}/\epsilon)
\,.
$$
Since this holds for all $t \ge t_0$ and
$\epsilon$ was arbitrary, (iii) follows.

Now assume that (ii) does not hold.
Let $\ev A_\ell$ be the event that the speed is at
least $\ell$, and note that $\ev A_\ell$ is $\omega$
measurable, by \ref l.nonrandom/.
Then by the famous bound of Varopoulos (1985) and Carne (1985), we have
$$
\lim_{t \to\infty} -{1 \over t} \log \mu^\omega_t\bigl(X(t)\bigr) 
\ge
\ell^2/2
$$
on $\ev A_\ell$.
In other words, (iii) does not hold.
\Qed

\procl c.nuqHD
\assumptions.
Suppose that $G$ is nonamenable.
Let $\omega$ be any $\Gamma$-invariant, ergodic,
insertion tolerant 
percolation
that has more than one infinite component a.s.
Then every infinite component of $\omega$ admits nonconstant bounded
harmonic Dirichlet functions.
\endprocl

\proof
By %
\ref t.cerg/ and ergodicity, it suffices to
establish the existence of nonconstant bounded harmonic Dirichlet
functions on $K(o)$ with positive probability.

We know from insertion tolerance and ergodicity that there
are infinitely many infinite components a.s. By 
insertion tolerance again, we also have that some, hence all, infinite
components have at least three ends. By \ref t.speed/, all infinite
components are transient. By insertion tolerance,
it follows that with positive probability,
$K(o)$ has a finite subset $K$
whose removal breaks $K(o)$ into at least two transient components.
In such a case, $K(o)$ has
nonconstant bounded harmonic Dirichlet functions
(e.g., the probability that a simple random walk
starting at $v$ eventually stays in a fixed transient
component of $K(o)\setminu K$ is such, as a function of $v$).
See Soardi (1994), Theorems 4.20 and 3.73.
This establishes our goal. \Qed

In order to prove transience in certain amenable cases,
we shall use:

\procl l.transtree
\procname{Transience of Big Trees}
If $T$ is any locally finite 
tree with $p_c(T) < 1$, then simple random walk is transient
on $T$.
\endprocl

\proof 
By Lyons (1990), the branching number of $T$ is $1/p_c(T)$ and this is the
critical value for transience of biased random walk on $T$. Since this is
larger than 1, it follows that, in particular, simple random walk is
transient. \Qed

For a Cayley graph $G$, let $\zeta_n$ be the number of elements
of $G$ at distance $n$ from $o$. It is evident that $\Seq{\zeta_n}$ is
submultiplicative, whence the
{\bf growth rate} $\gr(G) := \lim \zeta_n^{1/n} = \inf \zeta_n^{1/n}$ exists.

\procl t.gtra
\procname{Transience Above the Reciprocal Growth Rate} 
Let $G$ be a Cayley graph with $\gr(G)>1$, and let $p \in (1/\gr(G),1)$.
Then simple random walk is
transient on every infinite cluster of $p$-Bernoulli percolation a.s.
\endprocl

\proof
Let $\omega$ be $p$-Bernoulli percolation.
By ergodicity and indistinguishability of components, it suffices to prove
transience of $K(o)$ with positive probability. As shown in
Lyons (1995), there is a tree $T\subset G$ with $p_c(T) = 1/\gr(G)$.
This means that the component $\omega'$ of $o$ in $\omega\cap T$
has $p_c(\omega') < 1$
with positive probability, whence by \ref l.transtree/, $\omega$ is
transient with positive probability. By Rayleigh monotonicity, the same is
true of the component of $o$ in $G$. \Qed

The following conjecture would imply \ref g.psp/ (by taking $\omega := G$
and using \ref l.SEBperc/):

\procl g.ent
\procname{Monotonicity of Entropy}
\assumptions.
Let $\omega$ and $\omega'$ be two $\Gamma$-invariant percolations on $G$
such that $\omega' \subseteq \omega$. Then the asymptotic entropy of
delayed simple random walk on $\omega'$ is at most the asymptotic entropy
of delayed simple random walk on $\omega$.
\endprocl

The following conjecture for finite graphs can be shown
to imply \ref g.ent/.

\procl g.entfin
Let $G$ be a finite graph and $C : \edges(G) \to \R^+$.
Consider the continuous-time (reversible) Markov chain $\Seq{X(t)}$ on
$\verts(G)$ whose transition rate from $u$ to $v$ is $C(u, v)$.
Let $h_t(v, C)$ be the entropy of $X(t)$ when $X(0) = v$ and $h_t(C) :=
\sum_{v \in \verts(G)} h_t(v, C)$. 
Then for all $t$, given two functions $C'$ and $C$ with
$C'(e) \le C(e)$ for all $e \in \edges(G)$, we have $h_t(C') \le h_t(C)$.
\endprocl

Here is an equivalent formulation of this conjecture.
Given a matrix $B$, let $H(B)$ be the sum
of $-b_{i,j}\log b_{i,j}$ over all entries $b_{i,j}$ of the matrix.
Let ${\cal A}_n$ be the space of $n\times n$ real
symmetric matrices with non-negative off-diagonal
terms and with each row summing to zero.
Then a reformulation of \ref g.entfin/ is that
$H(\exp A)$ is (weakly) monotone increasing in the
off-diagonal entries of $A$, where $A$ ranges in ${\cal A}_n$.

\bsection{Harmonic Dirichlet Functions}{s.HD}

In this section, we study the existence of nonconstant harmonic Dirichlet
functions on percolation components.

We first describe the spanning forest measures we use.
A {\bf spanning tree} of a finite graph is a subgraph without
cycles that is connected and includes every vertex of the graph.  Motivated
by some questions of R.~Lyons, Pemantle (1991) showed that if an infinite
graph $G$ is exhausted by finite subgraphs $G_n$, then the uniform
distributions on the spanning trees of $G_n$ converge weakly to a measure
supported on spanning forests\ftnote{*}{In graph theory, ``spanning
forest'' usually means a maximal subgraph without cycles, i.e., a spanning
tree in each connected component. We mean, instead, a subgraph without
cycles that contains every vertex.} of $G$.  We call this the {\bf free
uniform spanning forest} ($\fusf$), since there is another natural
construction where the exterior of $G_n$ is identified to a single vertex
(``wired'') before passing to the limit.  This second construction, which
we call the {\bf wired uniform spanning forest} ($\wusf$), was implicit in
Pemantle's paper and was made explicit by H\"aggstr\"om (1995). Both
measures are concentrated on the set of forests, all of whose trees are
infinite. See
\BLPSusf\ or Lyons (1998) for an exposition and more details. For
convenience, we will use the symbols $\fusf$ and $\wusf$ also for the
uniform measure on spanning trees of a finite graph. For the proof of \ref
t.notOHD/, we shall have need of one more measure on infinite
transient graphs, the {\bf oriented wired
spanning forest}, denoted $\owsf$. We refer to \BLPSusf\ for its
definition. For our purposes, it is enough to know that it is the same as
$\wsf$, except that each edge in the forest is oriented in such a way that
there is exactly one outgoing edge from each vertex. All these measures,
$\fsf$, $\wsf$, and $\owsf$, are invariant under $\Aut(G)$.
We typically denote the random spanning forest by $\fo$.

The following three lemmas are taken from \BLPSusf:

\procl l.uniqusf
\procname{$\OHD$ Criterion}
For any (connected) graph $G$, we have $\fusf = \wusf$ iff $G \in \OHD$.
\endprocl

Write $\EGw$, $\EGf$ for expectation with respect to the random spanning
forests on $G$.

\procl l.equaldeg  %
\procname{Domination}
For any graph $G$, we have
$\EGf[\deg_\fo v ] \ge \EGw[\deg_\fo v ]$ for every
$v\in \vertex$, with equality for every $v$ iff $\fsf = \wsf$.
\endprocl

\procl l.ExpDeg2 %
\procname{$\wsf$-Expected Degree}
In any infinite transitive graph $G$, the 
$\wsf$-expected degree of every vertex is 2.
\endprocl

The following lemma is from \BLPSgip:

\procl l.degeq \procname{Small Trees and Expected Degree}
Let $\Gamma$ be a closed unimodular subgroup of
$\Aut(G)$ that acts transitively on $G$ and let $\omega$ be the
configuration of a
$\Gamma$-invariant percolation on $G$.  Fix a vertex $o$.  Let
$F_o$ be the event that $K(o)$ is an infinite tree with
finitely many ends, and let $F'_o$ be the event
that $K(o)$ is a finite tree.
\beginitems
\itemrm{(i)} If $P[F_o]>0$, then $\Ebig{D(o)\bigm|F_o}=2$.
\itemrm{(ii)} If $P[F'_o]>0$, then $\Ebig{D(o)\bigm|F'_o}<2$.
\enditems
\endprocl

Fix any basepoint $o \in \verts(G)$ and
let $ \ev A_o$ be the event that $K(o)$ is infinite.
Let $\E$ refer to the probability measure of the percolation.
Extending the notation above, we write $\Eof$ and $\Eow$
for expectation with respect to the free and
wired spanning forest measures on $\omega$ (given $\omega$).

\procl t.amen
\procname{$\OHD$ Stability when Amenable}
Let $G$ be an amenable graph with a transitive automorphism group
$\Gamma\subset\Aut(G)$ and $\omega$ a
$\Gamma$-invariant percolation.
Then a.s.\ every component of $\omega$ is in $\OHD$.
\endprocl

\proof We show that the measures $\fsf$ and $\wsf$ coincide on ${K(o)}$
a.s.\ given $ \ev A_o$.
By the argument in
Burton and Keane (1989), a.s.\ no component
in any invariant percolation on $G$ can
have more than 2 ends. Applying this to the percolations given by
taking the $\fsf$ or the $\wsf$ of each component (independently for each
component) of $\omega$ in
conjunction with \ref l.degeq/, we obtain that
$$
\Ebig{\Eow[\deg_\fo o] \bigm|  \ev A_o} 
= 2 =
\Ebig{\Eof[\deg_\fo o] \bigm|  \ev A_o} 
\,.
$$
Hence, the result follows from \ref l.equaldeg/.
\Qed

\procl l.recTree
\procname{Expected Degree for Recurrent Trees}
\assumptions.
Let $\omega$ be a $\Gamma$-invariant random forest
in $G$.  Suppose that a.s.\ all components of $\omega$ are recurrent.
Then $\E[\deg_\omega o]\leq 2$.
\endprocl

\proof This follows from Lemmas \briefref l.ends/,
\briefref l.degeq/, and \briefref l.transtree/. \Qed

\procl t.notOHD
\procname{$\neg \OHD$ Stability when High Marginals}
\assumptions.
If $G \notin \OHD$, then there is some $p_0 < 1$ such that for every
$\Gamma$-invariant bond percolation $\omega$ with $\inf_{e \in \edge} \P[e \in
\omega] > p_0$, some component of $\omega$ is not in $\OHD$ with positive
probability. If $\omega$ is ergodic and has indistinguishable components,
then the same hypotheses imply that a.s., no infinite component of
$\omega$ is in $\OHD$.
\endprocl

\proof We show that $p_0 := 2/\EGf[\deg_\fo o]$ works.

First, $p_0 < 1$ by Lemmas \briefref l.uniqusf/,
\briefref l.equaldeg/, and \briefref l.ExpDeg2/.
Let $\ev T_o$ be the event that $K(o)$, the component of
$o$ in $\omega$, is transient.
We claim that
$$
\P[\ev T_o] > 0
\label e.T_oHappens
$$
and
$$
\Ebig{\Eow[\deg_\fo o] \bigm|  \ev T_o} 
= 2 <
\Ebig{\Eof[\deg_\fo o] \bigm|  \ev T_o} 
\,.
\label e.toshow
$$
This suffices for the first statement
by \ref l.equaldeg/. The second statement then follows by ergodicity.

Let $\fo'$ be the union of components of $\fo$ that
are contained in recurrent components of $\omega$.
\ref l.recTree/ implies that if $\P[\ev T_o] < 1$, then
$\Eof[\deg_\fo o \mid o\in\fo'] \leq 2$,
which means $\Eof[\deg_\fo o \mid \neg  \ev T_o] \leq 2$.
Consequently, \ref e.T_oHappens/ and the inequality in
\ref e.toshow/ will be established once we prove
$$
\Ebig{\Eof[\deg_\fo o]}
>
p_0 \EGf[\deg_\fo o]
= 2\,.
$$
Now
$$
\Eof[\deg_\fo o]
=
\sum_{x \sim o} \Pof\big[ [o, x] \in \fo\big]
=
\sum_{x \sim o \atop [o, x] \in \omega} \Pof\big[ [o, x] \in \fo\big]
\,.
\label e.1
$$
Let $B_n$ be the ball of radius $n$ centered at $o$ in $G$.
By Kirchhoff's theorem and Rayleigh's monotonicity principle (see Lyons and
Peres (1998) or \BLPSusf), for
each $n$, and each $e\in\omega$,
$$
\P^{\omega \cap B_n}_\fsf[e \in \fo]
\geq
\P^{B_n}_\fsf[e \in \fo]
\,.
$$
Taking a limit as $n\to\infty$, we obtain
$$
\P^{\omega}_\fsf[e \in \fo]
\geq
\P_\fsf[e \in \fo]
$$
by the definition of the $\fusf$,
whence \ref e.1/ gives
$$
\Eof[\deg_\fo o]
\ge
\sum_{x \sim o \atop [o, x] \in \omega} \PGf[e \in \fo]
\,.
$$
Taking expectation, we obtain
\begineqalno
\Ebig{\Eof[\deg_\fo o]}
&\ge
\Eleft{\sum_{x \sim o} \I{[o, x] \in \omega} \PGf\big[[o,x] \in \fo\big]}
\cr&=
\sum_{x \sim o} \Pbig{[o, x] \in \omega} \PGf\big[[o,x] \in \fo\big]
\cr&>
\sum_{x \sim o} p_0 \PGf\big[[o,x] \in \fo\big]
\cr&=
p_0 \EGf[\deg_{\fo}o]
\,,
\endeqalno
as desired. 

For the equality in \ref e.toshow/, we use the oriented wired spanning
forest, $\owsf$, on each transient component of $\omega$, chosen
independently on each component.
Let $\varphi(x, y)$ be the probability that ($K(x)$ is transient and that)
$[x, y]$
belongs to the oriented wired spanning forest of $\omega$. Since $\owsf$ is
$\Gamma$-invariant, $\varphi$ is invariant under the diagonal action of
$\Gamma$, whence the Mass-Transport Principle says that
$$
\sum_x \varphi(o, x)
=
\sum_x \varphi(x, o)
\,.
$$
The left-hand side is the expected outdegree of $o$, which is $\P[\ev
T_o]$. Hence, the right-hand side, the expected in-degree of $o$, is also
$\P[\ev T_o]$. This shows that $\Ebig{\Eow[\deg_\fo o] \bigm|  \ev T_o} =
2$.
\Qed

\procl x.infinite
It does not suffice in \ref t.notOHD/ to assume merely that the components
of $\omega$ are infinite. For example, if $\omega$ is given by the $\wusf$,
then every component is a tree with one end (\BLPSusf), whence is recurrent
and in $\OHD$. However, as stated in \ref g.nhd/, we believe that this
{\it is\/} sufficient for Bernoulli percolation.
\endprocl

\procl x.nonunimod
The hypothesis that $\Gamma$ be unimodular cannot be omitted in \ref t.notOHD/.
For example, let $G$ be a regular tree of degree 3 and $\xi$ be an end of
$G$. Let $\Gamma$ be the group of automorphisms of $G$ that fix $\xi$.
Let $H_n$, $n\in\Z$, be the horocycles with respect to $\xi$.
(More precisely, fix a basepoint $o\in\verts$, let $\Seq{v_m}$ be a
sequence converging to $\xi$. Then a vertex $v$ is in $H_n$ iff
$\dist(v_m,v)-\dist(v_m,o)=n$ for all but finitely many $m$.)
To define $\omega$, we first define a percolation $\eta$.
Given any $p_0 < 1$, for each $n$ independently, let all the edges joining
$H_n$ to $H_{n+1}$ be in $\eta$ with probability $p_0$.
Each component of $\eta$ is a finite tree a.s. 
For each component $K$ of $\eta$, let $n(K)$ be the largest $n$ such that
$K \cap H_{n} \ne \emptyset$. Choose an edge joining $K \cap H_{n(K)}$
to $H_{n(K)+1}$ at random uniformly among all such
edges and independently for each $K$;
let $\eta'$ be the set of the chosen edges (over all $K$).
Now let $\omega := \eta \cup \eta'$. 
Each component of $\omega$ is a tree with exactly one end, so is recurrent
and in $\OHD$.
Yet $G \notin \OHD$. 
\endprocl

\procl q.nonuniber
Does \ref t.notOHD/ hold for Bernoulli percolation when the
unimodularity assumption is omitted?
\endprocl

\bsection{Anchored Expansion and Stability}{s.anchor}

Cheeger's inequality relates the isoperimetric constant, which is  
geometric, to the spectral radius, which governs the exponential decay of
return probabilities of simple random walk or
Brownian motion. We introduce
a geometric constant that we hope can replace the
isoperimetric constant in graphs and manifolds that are not uniformly
expanding.

Consider, for example, the hyperbolic space $\H^n$ and perturb the metric 
on an extremely
sparse sequence of balls with radii growing very slowly to infinity; for 
instance, pick the center of the $n$-th ball at distance $e^n$ from a fixed
origin, and let $\log n$ be its radius. If we modify the metric inside these
balls so that it is flat on sub-balls of half the radius, we get a manifold
with zero isoperimetric constant. Many
properties of $M$ (such as the existence of nonconstant bounded harmonic
functions or the speed of Brownian motion), are, however, unchanged
from $\H^n$.

\procl d.anchor 
Fix some basepoint $o \in G$. Call
$$
\anc(G) := \lim_{n \to \infty}\,
\inf \left\{{| \bde S| \over |S|} \st
o\in S\subset\verts(G),\,S \hbox{ is connected},\,n\leq |S|<\infty \right\}
$$   
the {\bf anchored expansion constant} of $G$.
Note that $\anc(G)$ is independent of the choice of the basepoint $o$
and that $\anc(G) \geq \isoe(G)$. 
\endprocl

By attaching a sequence of paths of length $1,2,\dots$
at a very sparse sequence of vertices of a binary tree,
we get an example of a graph $G$ for which $\anc(G) > \isoe(G)= 0$.

When the isoperimetric constant $\isoe(G)$ of a bounded
degree (not necessarily transitive) graph is positive,
Dodziuk's (1984) discrete version of Cheeger's inequality 
gives an upper bound $\overline{\rho}<1$ for the spectral radius $\rho(G)$,
where $\overline{\rho}$ depends on $\isoe(G)$ and the maximum degree in $G$.
In turn, this implies that the $\liminf$ speed
$\Lambda^-$ of simple random walk $X$ starting at a basepoint
$o\in\verts(G)$ is positive almost surely.

By analogy, we make the following

\procl g.anchorspeed Let $G$ be a bounded degree
graph with $\anc(G)>0$.  Then
$\Lambda^->0$ with positive probability.
\endprocl

It might even be the case that $\anc(G)>0$ implies
$\Lambda^->0$ a.s.

Thomassen (1992)
has shown that if a graph satisfies a certain ``rooted'' (= ``anchored'')
isoperimetric inequality, then it is transient.
\ref g.anchorspeed/ has a stronger hypothesis and a stronger conclusion.

The motivation for looking at $\anc(G)$ is that $\anc$ is
more stable than $\isoe$ under random perturbations of $G$.
For example, let $G$ be an infinite graph of bounded degree 
and pick a probability distribution $\P$ on 
the strictly positive integers. Replace each edge $e \in G$  by a path of 
length $L_e$, where $L_e$ is distributed according to $\P$, and
all $L_e$ ($e\in\edges(G)$) are independent.
Let $G^\P$ denote the random graph obtained in this way.
If $\P$ has a bounded support, then $\isoe(G)>0$ implies
$\isoe\left(G^\P\right)>0$,
while if $\P$ has unbounded support then, almost surely,
$\isoe\left(G^\P\right)=0$.

\procl q.hPa
Does $\isoe(G)>0$ imply that $\anc\left(G^\P\right)>0$ a.s.\ 
when $\P$ is the geometric distribution on the positive integers?
What about other distributions $\P$ with finite mean?
Can this be settled in the case where $G$ is a regular tree?
\endprocl

It is not hard to construct $\P$ and $G$ such that $\isoe(G)>0$
while $\anc\left(G^\P\right)=0$ a.s. (Take $G$ to be a binary tree and $\P$
to have a fat tail.)

Lyons, Pemantle, and Peres (1995) proved %
that simple random walk on a random perturbation of any
regular tree, with $\P$ the geometric distribution, has positive speed
almost surely.
More generally, simple random walk has positive speed on every
supercritical Galton-Watson tree a.s.\ given nonextinction.

\procl q.GW 
Is $\anc(T) > 0$ a.s.\ for supercritical Galton-Watson trees given
nonextinction?
\endprocl

\procl q.ancBern If $\isoe(G) > 0$ and $\omega$ is Bernoulli percolation on
$G$, must every infinite component $K$ of $\omega$ have $\anc(K)>0$ a.s.?
\endprocl

\beginreferences

Babson, E. \and Benjamini, I. (1998) Cut sets and normed cohomology with
applications to percolation, 
{\it Proc. Amer. Math. Soc.}, to appear.

Benjamini, I., Lyons, R., Peres, Y., \and Schramm, O. (1997)
Group-invariant percolation on graphs, {\it preprint}.

Benjamini, I., Lyons, R., Peres, Y., \and Schramm, O. (1998)
Uniform spanning forests, {\it preprint}.

Benjamini, I., Pemantle, R. \and Peres, Y. (1998)
Unpredictable paths and percolation, {\it Ann. Probab.}, to appear.

Benjamini, I. \and Schramm, O. (1996a) Harmonic functions on planar
and almost planar graphs and manifolds, via circle packings, {\it
Invent. Math.} {\bf 126}, 565--587.

Benjamini, I. \and Schramm, O. (1996b) Percolation beyond $\Z^d$, many
questions and a few answers, {\it Electronic Commun. Probab.} {\bf 1},
71--82.

Benjamini, I. \and Schramm, O. (1997) 
Every graph with a positive Cheeger constant contains a tree with
a positive Cheeger constant, {\it Geom. Funct. Anal.}
{\bf 7}, 403--419.

Burton, R. M. \and Keane, M. (1989)
Density and uniqueness in percolation, {\it Commun. Math.
Phys.} {\bf 121}, %
501--505.

Carne, T. K. (1985) A transmutation formula for Markov chains, {\it
Bull. Sci. Math.} (2) {\bf 109}, 399--405.

De Masi, A., Ferrari, P. A., Goldstein, S., \and Wick, W. D. (1989)
An invariance principle for reversible Markov processes. Applications to 
random motions in random environments, {\it J. Statist. Phys.} {\bf 55},
787--855. 

Derriennic, Y. (1980) Quelques applications du th\'eor\`eme ergodique
sous-additif, {\it Ast\'er\-isque} {\bf 74}, 183--201.

Dodziuk, J.  (1984) Difference equations, isoperimetric inequality, and
transience of certain random walks, {\it Trans. Amer. Math. Soc.} {\bf
284}, 787--794.

Grimmett, G. R. (1989) {\it Percolation}. Springer, New York.

Grimmett, G. R., Kesten, H. \and Zhang, Y. (1993) Random walk on the 
infinite cluster of the percolation model, {\it Probab. Theory Relat.
Fields} {\bf 96}, 33--44. 

Grimmett, G. R. \and Newman, C. M. (1990) Percolation in $\infty + 1$
dimensions, in {\it Disorder in Physical Systems}, G. R. Grimmett and D. J. A.
Welsh (editors), pp.~219--240. Clarendon Press, Oxford.

H\"aggstr\"om, O. (1995) Random-cluster measures and uniform spanning
trees, {\it Stoch. Proc. Appl.} {\bf 59}, 267--275.

H\"aggstr\"om, O. (1997) Infinite clusters in dependent automorphism
invariant percolation on trees, {\it Ann. Probab.} {\bf 25}, 1423--1436.

H\"aggstr\"om, O. \and Peres, Y. (1997) Monotonicity of uniqueness for
percolation on Cayley graphs: all infinite clusters are born
simultaneously, {\it preprint}.

Holopainen, I. \and Soardi, P. M. (1997)
$p$-harmonic functions on graphs and manifolds,
{\it Manuscripta Math.} {\bf 94}, 95--110.

Kaimanovich, V. A. (1990) Boundary and entropy of random walks in random
environment, in {\it Probability theory and mathematical statistics}, Vol. I
(Vilnius, 1989), 573--579, ``Mokslas", Vilnius, 1990.

Kaimanovich, V. A. (1992) Dirichlet norms, capacities and generalized
isoperimetric inequalities for Mar\-kov operators, {\it Potential Anal.}
{\bf 1}, 61--82.

Kaimanovich, V. A. \and Vershik, A. M. (1983) Random walks on
discrete groups: boundary and entropy, {\it Ann. Probab.} {\bf 11},
457--490.

Kaimanovich, V. A. \and Woess, W. (1998) The Poisson boundary of
quasi-transitive graphs, {\it in preparation.}

Kanai, M. (1986) Rough isometries and the parabolicity of riemannian
manifolds, {\it J. Math. Soc. Japan} {\bf 38}, 227--238.

Kesten, H. (1959) Symmetric random walks on groups, {\it Trans. Amer. Math.
Soc.} {\bf 92}, 336--354.

Lyons, R. (1990) Random walks and percolation on trees, {\it Ann. Probab.}
{\bf 18}, 931--958.

Lyons, R. (1995) Random walks and the growth of groups, {\it C. R. Acad.
Sci. Paris} {\bf 320}, 1361--1366.

Lyons, R. (1998) A bird's-eye view of uniform spanning trees and forests,
in {\it Microsurveys in Discrete Probability}, D. Aldous and J. Propp
(editors). AMS, Providence, to appear.
 
Lyons, R., Pemantle, R. \and Peres, Y. (1995) Ergodic theory on
Galton-Watson trees: speed of random walk and dimension of harmonic
measure, {\it Ergodic Theory Dynamical Systems} {\bf 15}, 593--619.

Lyons, R. \and Peres, Y. (1998) {\it Probability on Trees and Networks}.
Cambridge University Press, in preparation. Current version available 
\ifsubmit
at {\tt http://php.indiana.edu/\~{}rdlyons/}.
\else
at \hfill\break {\tt http://php.indiana.edu/\~{}rdlyons/}.
\fi

Lyons, R. \and Schramm, O. (1998) {\it Paper in preparation}.

Lyons, T. (1987) Instability of the Liouville property for quasi-isometric
Riemannian manifolds and reversible Markov chains, {\it J. Diff. Geom.}
{\bf 26}, 33--66.

Medolla, G. \and Soardi, P. M. (1995) Extension of Foster's
averaging formula to infinite networks with moderate growth, {\it Math. Z.}
{\bf 219}, 171--185.

Pemantle, R. (1991) Choosing a spanning tree for the integer lattice
uniformly, {\it Ann. Probab.} {\bf 19}, 1559--1574.

\vbox{
Petersen, K. (1983) {\it Ergodic Theory.} Cambridge University Press,
Cambridge.
}

Saloff-Coste, L. (1992)
A note on Poincare, Sobolev, and Harnack inequalities,
{\it Internat. Math. Res. Notices}, no. 2, 27--38.

Soardi, P. M. (1993) Rough isometries and Dirichlet finite harmonic
functions on graphs, {\it Proc. Amer. Math. Soc.} {\bf 119}, 1239--1248.

Soardi, P. M. (1994) {\it Potential Theory on Infinite Networks.}
Springer, Berlin.  %

Soardi, P. M. \and Woess, W. (1990) Amenability,
unimodularity, and the spectral radius of random walks on infinite  
graphs, {\it Math. Z.} {\bf 205}, 471--486. 

Thomassen, C. (1992) Isoperimetric inequalities and transient random walks on
graphs, {\it Ann.\ Probab.{}} {\bf 20}, {1592--1600}.

Varopoulos, N. Th. (1985) Long range estimates for Markov chains, {\it Bull.
Sci. Math.} (2) {\bf 109}, 225--252.

Woess, W. (1994) Random walks on infinite graphs and groups --- a survey on
selected topics, {\it Bull. London Math. Soc.} {\bf 26}, 1--60.

\endreferences

\filbreak
\begingroup
\eightpoint\sc
\parindent=0pt\baselineskip=10pt

\def\emailwww#1#2{\par\qquad {\tt #1}\par\qquad {\tt #2}\smallskip}

Mathematics Department,
The Weizmann Institute of Science,
Rehovot 76100, Israel
\emailwww{itai@wisdom.weizmann.ac.il}
{http://www.wisdom.weizmann.ac.il/\~{}itai/}

Department of Mathematics,
Indiana University,
Bloomington, IN 47405-5701, USA
\emailwww{rdlyons@indiana.edu}
{http://php.indiana.edu/\~{}rdlyons/}

Mathematics Department,
The Weizmann Institute of Science,
Rehovot 76100, Israel
\emailwww{schramm@wisdom.weizmann.ac.il}
{http://www.wisdom.weizmann.ac.il/\~{}schramm/}

\endgroup

\bye